# SELFADJOINT EXTENSIONS OF A SINGULAR MULTİPOİNT DIFFERENTIAL OPERATOR FOR FİRST ORDER

by

## Z. I. ISMAILOV[1]  & R. ÖZTÜRK MERT[1]

[1] Karadeniz Technical University, Faculty of Sciences, Department of Mathematics, 61080, Trabzon, TURKEY; e-mail: zameddin@yahoo.com;

**Abstract**

In this work, firstly in the direct sum of Hilbert spaces of vector-functions $L^2(H,(-\infty,a_1)) \oplus L^2(H,(a_2,b_2)) \oplus L^2(H,(a_3,+\infty))$, $-\infty < a_1 < a_2 < b_2 < a_3 < +\infty$ all selfadjoint extensions of the minimal operator generated by linear singular multipoint symmetric differential expression $l = (l_1, l_2, l_3), l_k = i\frac{d}{dt} + A_k$ with a selfadjoint operator coefficient $A_k$ $k = 1,2,3$ in any Hilbert space $H$, are described in terms of boundary values. Later structure of the spectrum of these extensions is investigated.

## 1. Introduction

Many problems arising in the modelling of processes of multi-particle quantum mechanics, in the quantum field theory, in the multipoint boundary value problems for the differential equations, in the physics of rigid bodies and ets support to study selfadjoint extension of symmetric differential operators in direct sum of Hilbert spaces([1-3]).

The general theory of selfadjoint extensions of symmetric operators in any Hilbert space and their spectral theory deeply have been investigated by many mathematicians ([4-7]). Applications of this theory to two point differential operators in Hilbert space of functions are continied today even.

It is well-known that for the existence of selfadjoint extension of the any linear closed densely defined symmetric operator $B$ in a Hilbert space $\mathcal{H}$, the necessary and sufficient condition is a equality of deficiency indices $m(B) = n(B)$, where $m(B) = \dim ker(B^* + i)$, $n(B) = \dim ker(B^* - i)$.

But the multipoint situations may occur in different tables in the following sense. Let $B_1$ and $B_2$ be minimal operators generated by the linear differential expression $i\frac{d}{dt}$ in the Hilbert space of functions $L^2(-\infty, a)$ and $L^2(b, +\infty), a < b$, respectively. In this case it is known rhat
$$(m(B_1), n(B_1)) = (0,1), (m(B_2), n(B_2)) = (1,0).$$

Consequently, $B_1$ and $B_2$ are maximal symmetric operators, but are not a selfadjoint. However, direct sum $B = B_1 \oplus B_2$ of operators in a direct sum $\mathcal{H} = L^2(-\infty, a) \oplus L^2(b, +\infty)$ spaces have an equal defect numbers $(1,1)$. Then by the general theory [4] it has a selfadjoint extension. On the other hand it can be easily shown in the form that
$$u_2(b) = e^{i\varphi} u_1(a), \varphi \in [0, 2\pi), u = (u_1, u_2), u_1 \in D(B_1^*), u_2 \in D(B_2^*).$$
In singular cases although to date no investigation so far has not been. But the physical and technical processes, many of the problems resulting from examination of the solution is of great importance for the singular cases.

In this work in second section by the method of J.W.Calkin- M.L.Gorbachuk theory all selfadjoint extensions of the minimal operator generated by singular multipoint symmetric differential operator for first order in the direct sum of Hilbert space $L^2(H, (-\infty, a_1)) \oplus L^2(H, (a_2, b_2)) \oplus L^2(H, (a_3, +\infty))$, $-\infty < a_1 < a_2 < b_2 < a_3 < +\infty$ in terms of boundary values are described. In third section the spectrum of such extensions is researched.

## 2. Description of Selfadjoint Extensions

Let $H$ be a separable Hilbert space and $a_1, a_2, b_2, a_3 \in R, a_1 < a_2 < b_2 < a_3$. In the Hilbert space $L^2(H, (-\infty, a_1)) \oplus L^2(H, (a_2, b_2)) \oplus L^2(H, (a_3, +\infty))$ of vector-functions consider the following linear multipoint differential expression in form
$$l(u) = (l_1(u_1), l_2(u_2), l_3(u_3)) = (iu_1' + A_1 u_1, iu_2' + A_2 u_2, iu_3' + A_3 u_3), \; u = (u_1, u_2, u_3),$$
where $A_k: D(A_k) \subset H \to H, k = 1,2,3$ is linear selfadjoint operators in $H$. In the linear manifold $D(A_k) \subset H$ introduce the inner product in form
$$(f, g)_{k,+} := (A_k f A_k, g)_H + (f, g)_H, \; f, g \in D(A_k), k = 1,2,3.$$
For the $k = 1,2,3$, $D(A_k)$ is a Hilbert space under the positive norm $\|\cdot\|_{k,+}$ respect to Hilbert space $H$. It is denoted by $H_{k,+}$. Denote the $H_{k,-}$ a Hilbert space with negative norm. It is clear that a operator $A_k$ is continuous from $H_{k,+}$ to $H$ and it's adjoint operator $\tilde{A}_k: H \to H_{k,-}$ is a extension of the operator $A_k$. On the other hand, the operator $\tilde{A}_k: D(\tilde{A}_k) = H \subset H_{k,-1} \to H_{k,-1}$ is a linear selfadjoint.

Now in the direct sum $L^2(H, (-\infty, a_1)) \oplus L^2(H, (a_2, b_2)) \oplus L^2(H, (a_3, +\infty))$ define by
$$\tilde{l}(u) = (\tilde{l}_1(u_1), \tilde{l}_2(u_2), \tilde{l}_3(u_3)), \tag{2.1}$$
where $u = (u_1, u_2, u_3)$ and $\tilde{l}_1(u_1) = iu_1' + \tilde{A}_1 u_1$, $\tilde{l}_2(u_2) = iu_2' + \tilde{A}_2 u_2, \tilde{l}_2(u_2) = iu_2' + \tilde{A}_2 u_2, \tilde{l}_3(u_3) = iu_3' + \tilde{A}_3 u_3$.

The minimal $L_{10}(L_{20} and L_{30})$ and maximal $L_1(L_2 and L_3)$ operators generated by differential expression $\tilde{l}_1(\tilde{l}_2 and \tilde{l}_3)$ in $L^2(H, (-\infty, a_1))(L^2(H, (a_2, b_2)) and L^2(H, (b, +\infty)))$ have been investigated in [8].

The operators defined by $L_0 = L_{10} \oplus L_{20} \oplus L_{30}$ and $L = L_1 \oplus L_2 \oplus L_3$ in the space $L^2 = L^2(H, (-\infty, a_1)) \oplus L^2(H, (a_2, b_2)) \oplus L^2(H, (a_3, +\infty))$ are called minimal and maximal (multipoint) operators generated by the differential expression (2.1), respectively. Note that the operator $L_0$ is a symmetric and $L_0^* = L$ in $L^2$. On the other hand, it is clear that, $m(L_{10}) = 0, n(L_{10}) = dim H, m(L_{20}) = dim H, n(L_{20}) = dim H, m(L_{30}) = dim H, n(L_{30}) = 0$.

Consequently, $m(L_0) = n(L_0) = 2 dim H > 0$. So the minimal operator $L_0$ has a selfadjoint extension ([4]). For example, the differential expression $\tilde{l}(u)$ with a boundary condition $u(a_1) = u(a_3)$, $u(a_2) = u(b_2)$ generates a selfadjoint operator in $L^2$.

Here it is described all selfadjoint extensions of the minimal operator $L_0$ in $L^2$ in terms of the boundary values.

Note that a space of boundary values has an important role in the theory of selfadjoint extensions of the linear symmetric differential operators ( [6,7] ).

Let $B: D(B) \subset \mathcal{H} \to \mathcal{H}$ be a closed densely defined symmetric operator in the Hilbert space $\mathcal{H}$, having equal finite or infinite deficiency indices. A triplet $(\mathfrak{H}, \gamma_1, \gamma_2)$, where $\mathfrak{H}$ is a Hilbert space, $\gamma_1$ and $\gamma_2$ are linear mappings of $D(B^*)$ into $\mathfrak{H}$, is called a space of boundary values for the operator $B$ if for any $f, g \in D(B^*)$
$$(B^*f, g)_\mathcal{H} - (f, B^*g)_\mathcal{H} = (\gamma_1(f), \gamma_2(g))_\mathfrak{H} - (\gamma_2(f), \gamma_1(g))_\mathfrak{H},$$
while for any $F, G \in \mathfrak{H}$, there exists an element $f \in D(B^*)$, such that $\gamma_1(f) = F$ and $\gamma_2(f) = G$.

Note that any symmetric operator with equal deficiency indices have at least one space of boundary values ([6]).

In first note that the following proposition which validity of this claim can be easily proved.

**Lemma 2.1:** The triplet $(H, \gamma_1, \gamma_2)$,
$$\gamma_1: D((L_{10} \oplus 0 \oplus L_{30})^*) \to H, \qquad \gamma_1(u) = \frac{1}{i\sqrt{2}}(u_1(a_1) + u_3(a_3)),$$
$$\gamma_2: D((L_{10} \oplus 0 \oplus L_{30})^*) \to H, \qquad \gamma_2(u) = \frac{1}{\sqrt{2}}(u_1(a_1) - u_3(a_3)),$$
$$u = (u_1, u_2, u_3) \in D((L_{10} \oplus 0 \oplus L_{30})^*)$$
is a space of boundary values of the minimal operator $L_{10} \oplus 0 \oplus L_{30}$ in the direct sum
$$L^2(H, (-\infty, a_1)) \oplus 0 \oplus L^2(H, (a_3, +\infty))$$

**Proof:** For the arbitrary $u = (u_1, u_2, u_3)$ and $v = (v_1, v_2, v_3)$ from $D(((L_{10} \oplus 0 \oplus L_{30})^*))$ validity the following equality
$$(Lu, v)_{L^2(H,(-\infty,a_1))\oplus 0 \oplus L^2(H,(a_3,\infty))} - (u, Lv)_{L^2(H,(-\infty,a_1))\oplus 0 \oplus L^2(H,(a_3,\infty))}$$
$$= (\gamma_1(u), \gamma_2(v))_H - (\gamma_2(u), \gamma_1(v))_H$$
can be easily verified. Now give any elements $f, g \in H$. Find the function $u = (u_1, u_2, u_3) \in D(((L_{10} \oplus 0 \oplus L_{30})^*))$ such that
$$\gamma_1(u) = \frac{1}{i\sqrt{2}}(u_1(a_1) + u_3(a_3)) = f \text{ and } \gamma_2(u) = \frac{1}{\sqrt{2}}(u_1(a_1) - u_3(a_3)) = g$$
that is,
$$u_1(a_1) = (if + g)/\sqrt{2} \text{ and } u_3(a_3) = (if - g)/\sqrt{2}.$$
If choose these functions $u_1(t), u_3(t)$ in following form
$$u_1(t) = \int_{-\infty}^{t} e^{s-a_1} ds (if + g)/\sqrt{2}, \; t < a_1,$$
$$u_2(t) = 0, \qquad a_2 < t < b_2,$$
$$u_3(t) = \int_{t}^{\infty} e^{a_3 - t} ds (if - g)/\sqrt{2}, \; t > a_3,$$
then it is clear that $(u_1, u_2, u_3) \in D(((L_{10} \oplus 0 \oplus L_{30})^*))$ and $\gamma_1(u) = f$, $\gamma_2(u) = g$.
□ Furthermore, using the result which obtained in [8] it is proved the next assertion.

**Lemma 2.2.** The triplet $(H, \Gamma_1, \Gamma_2)$,
$$\Gamma_1: D((0 \oplus L_{20} \oplus 0)^*) \to H, \; \Gamma_2(u) = \frac{1}{i\sqrt{2}}(u_2(a_2) + u_2(b_2)),$$
$$\Gamma_2: D((0 \oplus L_{20} \oplus 0)^*) \to H, \Gamma_2(u) = \frac{1}{\sqrt{2}}(u_2(a_2) - u_2(b_2)),$$
$$u = (u_1, u_2, u_3) \in D((0 \oplus L_{20} \oplus 0)^*)$$

is a space of boundary values of the minimal operator $0\oplus L_0\oplus 0$ in the direct sum $0\oplus L^2(H,(a_2,b_2))\oplus 0$.

On the other hand can be easily established the following result.

**Lemma 2.3.** Every selfadjoint extension of $L_0$ in $L^2 = L^2(H,(-\infty,a_1))\oplus L^2(H,(a_2,b_2))\oplus L^2(H,(a_3,+\infty))$ is a direct sum of selfadjoint extensions of the minimal operator $L_{10}\oplus 0\oplus L_{30}$ in $L^2(H,(-\infty,a_1))\oplus 0\oplus L^2(H,(a_3,+\infty))$ and minimal operator $0\oplus L_0\oplus 0$ in $0\oplus L^2(H,(a_2,b_2))\oplus 0$.

Finally, using the method in [6] can be established the following result.

**Theorem 2.4:** If $\tilde{L}$ is a selfadjoint extension of the minimal operator $L_0$ in $L^2$, then it generates by the differential expression (2.1) and the boundary conditions
$$u_3(a_3) = W_1 u_1(a_1),$$
$$u_2(b_2) = W_2 u_2(a_2),$$
where $W_1, W_2: H \to H$ are a unitary operators. Moreover, the unitary operators $W_1, W_2$ in $H$ are determined uniquely by the extension $\tilde{L}$, i.e. $\tilde{L} = L_{W_1 W_2}$ and vice versa.

## 3. The Spectrum of the Selfadjoint Extensions

In this section the structure of the spectrum of the selfadjoint extension $L_{W_1 W_2}$ in $L^2$ will be investigated.

In this case by the Lemma 2.3. it is clear that
$$L_{W_1 W_2} = L_{W_1}\oplus L_{W_2},$$
where $L_{W_1}$ and $L_{W_2}$ are selfadjoint extensions of the minimal operators $L_0(1,0,1) = L_{10}\oplus 0\oplus L_{30}$ and $L_0(0,1,0) = 0\oplus L_0\oplus 0$ in the Hilbert spaces $L^2(1,0,1) = L^2(H,(-\infty,a_1))\oplus 0\oplus L^2(H,(a_3,+\infty))$ and $L^2(0,1,0) = 0\oplus L^2(H,(a_2,b_2))\oplus 0$ respectively.

First of all, we have to prove the following result.

**Theorem 3.1:** The point spectrum of any selfadjoint extension $L_{W_1}$, in the Hilbert space $L^2(1,0,1)$ is empty, i.e.,
$$\sigma_p(L_{W_1}) = \emptyset.$$

**Proof:** In first consider the following problem for the spectrum of the selfadjoint extension $L_{W_1}$ of the minimal operator $L_0(1,0,1)$ in the Hilbert space $L^2(1,0,1)$
$$L_{W_1} u = \lambda u, u = (u_1, 0, u_3) \in L^2(1,0,1),$$
that is,
$$\tilde{l}_1(u_1) = iu_1' + \tilde{A}_1 u_1 = \lambda u_1, \quad u_1 \in L^2(H,(-\infty,a_1)),$$
$$\tilde{l}_3(u_3) = iu_3' + \tilde{A}_3 u_3 = \lambda u_3, \quad u_3 \in L^2(H,(a_3,+\infty)), \quad \lambda \in R$$
$$u_3(a_3) = W_1 u_1(a_1) \quad .$$

The general solution of this problem is
$$u_1(t) = e^{i(\tilde{A}_1-\lambda)(t-a_1)} f_1^*, \quad t < a_1$$
$$u_3(t) = e^{i(\tilde{A}_3-\lambda)(t-a_3)} f_3^*, \quad t > a_3$$
$$f_3^* = W_1 f_1^*, \quad f_1^*, f_3^* \in H.$$

It is clear that for the $f_1^* \neq 0, f_3^* \neq 0$ the functions $u_1 \notin L^2(H,(-\infty,a_1))$, $u_2 \notin L^2(H,(a_3,\infty))$. So for every unitary operator $W_1$ we have $\sigma_p(L_{W_1}) = \emptyset$.

Since residual spectrum of any selfadjoint operator in any Hilbert space is empty, then furthermore the continuous spectrum of the selfadjoint extensions $L_{W_1}$ of the minimal operator $L_0(1,0,1)$ in the Hilbert space $L^2(1,0,1)$ is investigated.

**Theorem 3.2:** Continuous spectrum of the any selfadjoint extension $L_{W_1}$ of the minimal operator $L_0(1,0,1)$ in the Hilbert space $L^2(1,0,1)$ is in form
$$\sigma_c(L_{W_1}) = \mathbb{R}.$$

**Proof.** Now before of all it will be researched the resolvent operator of the extension $L_{W_1}$ generated by the differential expression $(\tilde{l}_1, 0, \tilde{l}_3)$ and the boundary condition
$$u_3(a_3) = W_1 u_1(a_1)$$
in the Hilbert space $L^2(1,0,1)$, i.e.

$$\tilde{l}_1(u_1) = iu'_1 + \tilde{A}_1 u_1 = \lambda u_1 + f_1, \quad u_1, f_1 \in L^2(H, (-\infty, a_1)),$$
$$\tilde{l}_3(u_3) = iu'_3 + \tilde{A}_3 u_3 = \lambda u_3 + f_3, \quad u_3, f_3 \in L^2(H, (a_3, +\infty)), \quad \lambda \in \mathbb{C}, \lambda_i = Im\lambda > 0$$
$$u_3(a_3) = W_1 u_1(a_1) \tag{3.1}$$

Now we will be shown that the following function
$$u(\lambda; t) = (u_1(\lambda; t), 0, u_3(\lambda; t)),$$
where
$$u_1(\lambda; t) = e^{i(\tilde{A}_1 - \lambda)(t - a_1)} f_1^* + i \int_t^{a_1} e^{i(\tilde{A}_1 - \lambda)(t-s)} f_1(s) ds, \quad t < a_1,$$
$$u_3(\lambda; t) = i \int_t^{\infty} e^{i(\tilde{A}_3 - \lambda)(t-s)} f_3(s) ds, \quad t > a_3,$$
$$f_1^* = W^* \left( i \int_{a_3}^{\infty} e^{i(\tilde{A}_3 - \lambda)(t-s)(b-s)} f_3(s) ds \right)$$

is a solution of the boundary value problems (3.1) in the Hilbert space $L^2(1,0,1)$. For this, it is sufficient to show that
$$u_1(\lambda; t) \in L^2(H, (-\infty, a_1)),$$
$$u_3(\lambda; t) \in L^2(H, (a_3, +\infty))$$
for the case $\lambda_i > 0$. Indeed, in this case

$$\|f_1^*\|_H^2 = \left\| \int_{a_3}^{\infty} e^{i(\tilde{A}_3 - \lambda)(a_3 - s)} f(s) ds \right\|_H^2 \leq \left( \int_{a_3}^{\infty} e^{\lambda_i(a_3 - s)} \|f(s)\|_H ds \right)^2$$
$$\leq \left( \int_{a_3}^{\infty} e^{2\lambda_i(a_3 - s)} ds \right) \left( \int_{a_3}^{\infty} \|f(s)\|_H^2 ds \right) = \frac{1}{2\lambda_i} \|f\|_{L^2(H,(a_3,+\infty))}^2 < \infty,$$

$$\left\| e^{i(\tilde{A}_1 - \lambda)(t - a_1)} f_1^* \right\|_{L^2(H,(-\infty,a_1))}^2 = \left\| e^{-i\lambda(t - a_1)} f_1^* \right\|_{L^2(H,(-\infty,a_1))}^2 = \int_{-\infty}^{a_1} \left\| e^{-i\lambda(t - a_1)} f_1^* \right\|_H^2 dt$$
$$= \int_{-\infty}^{a_1} e^{2\lambda_i(t - a_1)} dt \|f_1^*\|_H^2 = \frac{1}{2\lambda_i} \|f_1^*\|_H^2 < \infty$$

and

$$\left\|i\int_t^{a_1} e^{i(\tilde{A}_1-\lambda)(t-s)}f_1(s)ds\right\|^2_{L^2(H,(-\infty,a_1))} \leq \int_{-\infty}^{a_1}\left(\int_t^{a_1} e^{\lambda_i(t-s)}\|f_1(s)\|_H ds\right)^2 dt$$

$$\leq \int_{-\infty}^{a_1}\left(\int_t^{a_1} e^{\lambda_i(t-s)}ds\right)\left(\int_t^{a_1} e^{\lambda_i(t-s)}\|f_1(s)\|^2 ds\right) dt$$

$$= \frac{1}{\lambda_i}\int_{-\infty}^{a_1}\int_t^{a_1} e^{\lambda_i(t-s)}\|f_1(s)\|^2 ds\, dt = \frac{1}{\lambda_i}\int_{-\infty}^{a_1}\left(\int_{-\infty}^{s} e^{\lambda_i(t-s)}\|f_1(s)\|^2 dt\right)ds$$

$$= \frac{1}{\lambda_i}\int_{-\infty}^{a_1}\left(\int_{-\infty}^{s} e^{\lambda_i(t-s)}dt\right)\|f_1(s)\|^2 ds = \frac{1}{\lambda_i^2}\int_{-\infty}^{a_1}\|f_1(s)\|^2 ds$$

$$= \frac{1}{\lambda_i^2}\|f_1\|^2_{L^2(H,(-\infty,a_1))} < \infty.$$

Furthermore,

$$\left\|i\int_t^{\infty} e^{i(\tilde{A}_3-\lambda)(t-s)}f_3(s)ds\right\|^2_{L^2(H,(a_3,+\infty))} \leq \int_{a_3}^{\infty}\left(\int_t^{\infty} e^{\lambda_i(t-s)}\|f_3(s)\|_H ds\right)^2 dt$$

$$\leq \int_{a_3}^{\infty}\left(\int_t^{\infty} e^{\lambda_i(t-s)}ds\right)\left(\int_t^{\infty} e^{\lambda_i(t-s)}\|f_3(s)\|^2 ds\right) dt$$

$$= \frac{1}{\lambda_i}\int_{a_3}^{\infty}\left(\int_t^{\infty} e^{\lambda_i(t-s)}\|f_3(s)\|^2 ds\right) dt = \frac{1}{\lambda_i}\int_{a_3}^{\infty}\left(\int_{a_3}^{s} e^{\lambda_i(t-s)}\|f_3(s)\|^2 dt\right) ds$$

$$= \frac{1}{\lambda_i}\int_{a_3}^{\infty}\left(\int_{a_3}^{s} e^{\lambda_i(t-s)}dt\right)\|f_3(s)\|^2 ds = \frac{1}{\lambda_i^2}\left(\int_{a_3}^{\infty}(1-e^{\lambda_i(a_3-s)})\|f_3(s)\|^2 ds\right)$$

$$\leq \frac{1}{\lambda_i^2}\|f_3\|^2_{L^2(H,(a_3,+\infty))} < \infty.$$

From above calculations imply that $u_1(\lambda;t) \in L^2(H,(-\infty,a_1))$, $u_3(\lambda;t) \in L^2(H,(a_3,+\infty))$ for $\lambda \in \mathbb{C}$, $\lambda_i = Im\lambda > 0$. On the other hand it can be to easy to verify that $u(\lambda;t) = (u_1(\lambda;t), 0, u_3(\lambda;t))$ is a solution of the boundary value problem (3.1).

In the case when $\lambda \in \mathbb{C}$, $\lambda_i = Im\lambda < 0$ solution of the boundary value problem
$$L_{W_1}u = \lambda u + f, \qquad u = (u_1,0,u_3), \qquad f = (f_1,0,f_3) \in L^2(1,0,1)$$
$$u_3(a_3) = W_1 u_1(a_1),$$
where $W_1$ is a unitary operator in $H$, is in the form $u(\lambda;t) = (u_1(\lambda;t), 0, u_3(\lambda;t))$,

$$\begin{cases} u_1(\lambda;t) = -i\int_{-\infty}^{t} e^{i(\tilde{A}_1-\lambda)(t-s)}f_1(s)ds, & t < a_1 \\ \\ u_3(\lambda;t) = e^{i(\tilde{A}_3-\lambda)(t-a_3)}f_3^* - i\int_{a_3}^{t} e^{i(\tilde{A}_3-\lambda)(t-s)}f_3(s)ds, & t > a_3, \end{cases}$$

$$f_3^* = W\left(-i\int_{-\infty}^{a_1} e^{i(\tilde{A}_1-\lambda)(a_1-s)}f_1(s)ds\right).$$

In first prove that $u(\lambda;t) \in L^2(1,0,1)$. In this case

$$\|u_1(\lambda;t)\|^2_{L^2(H,(-\infty,a_1))} = \int_{-\infty}^{a_1} \left\| -i \int_{-\infty}^{t} e^{i(\tilde{A}_1-\lambda)(t-s)} f_1(s) ds \right\|^2_H dt$$

$$\leq \int_{-\infty}^{a_1} \left( \int_{-\infty}^{t} e^{\lambda_i(t-s)} ds \right) \left( \int_{-\infty}^{t} e^{\lambda_i(t-s)} \|f_1(s)\|^2_H ds \right) dt$$

$$= \frac{1}{|\lambda_i|} \int_{-\infty}^{a_3} \int_{-\infty}^{t} e^{\lambda_i(t-s)} \|f_1(s)\|^2_H \, ds \, dt$$

$$= \frac{1}{|\lambda_i|} \int_{-\infty}^{a_1} \left( \int_{s}^{a_1} e^{\lambda_i(t-s)} \|f_1(s)\|^2_H dt \right) ds = \frac{1}{|\lambda_i|} \int_{-\infty}^{a_1} \left( \int_{s}^{a_1} e^{\lambda_i(t-s)} dt \right) \|f_1(s)\|^2_H ds$$

$$= \frac{1}{|\lambda_i|^2} \int_{-\infty}^{a_1} (1 - e^{\lambda_i(a_1-s)}) \|f_1(s)\|^2_H \, ds \leq \frac{1}{|\lambda_i|^2} \|f_1\|^2_{L^2(H,(-\infty,a_1))} < \infty,$$

$$\|f_3^*\|^2_H = \left\| \int_{-\infty}^{a_1} e^{i(\tilde{A}_1-\lambda)(a_1-s)} f_1(s) ds \right\|^2_H \leq \left( \int_{-\infty}^{a_1} e^{\lambda_i(a_1-s)} \|f_1(s)\|_H ds \right)^2$$

$$\leq \left( \int_{-\infty}^{a_1} e^{2\lambda_i(a_1-s)} ds \right) \left( \int_{-\infty}^{a_1} \|f_1(s)\|^2_H \, ds \right)$$

$$= \frac{1}{2|\lambda_i|} \|f_1\|^2_{L^2(H,(-\infty,a_1))} < \infty,$$

$$\left\| e^{i(\tilde{A}_3-\lambda)(t-a_3)} f_3^* \right\|^2_{L^2(H,(a_3,+\infty))} \leq \int_{a_3}^{\infty} e^{2\lambda_i(t-a_3)} dt \, \|f_3^*\|^2_H = \frac{1}{2|\lambda_i|} \|f_3^*\|^2_H$$

$$\leq \frac{1}{4|\lambda_i|^2} \|f\|^2_{L^2(H,(a_3,+\infty))} < \infty$$

and

$$\left\| \int_{a_3}^{t} e^{i(\tilde{A}_3-\lambda)(t-s)} f_3(s) ds \right\|^2_{L^2(H,(a_3,+\infty))} \leq \int_{a_3}^{\infty} \left( \int_{a_3}^{t} e^{\lambda_i(t-s)} \|f_3(s)\|_H ds \right)^2 dt$$

$$\leq \int_{a_3}^{\infty} \left( \int_{a_3}^{t} e^{\lambda_i(t-s)} ds \right) \left( \left( \int_{a_3}^{t} e^{\lambda_i(t-s)} \|f_3(s)\|^2_H ds \right) \right) dt$$

$$= \int_{a_3}^{\infty} \left( \frac{1}{\lambda_i} (1 - e^{\lambda_i(t-a_3)}) \right) \left( \int_{a_3}^{t} e^{\lambda_i(t-s)} \|f_3(s)\|^2_H ds \right) dt$$

$$\leq \frac{1}{|\lambda_i|} \int_{a_3}^{\infty} \left( \int_{a_3}^{t} e^{\lambda_i(t-a_3)} \|f_3(s)\|^2_H ds \right) dt$$

$$= \frac{1}{|\lambda_i|} \int_{a_3}^{\infty} \left( \int_{s}^{\infty} e^{\lambda_i(t-s)} \|f_3(s)\|^2_H dt \right) ds$$

$$= \frac{1}{|\lambda_i|} \int_{a_3}^{\infty} \left( \int_s^{a_3} e^{\lambda_i(t-s)} dt \right) \|f_3(s)\|_H^2 ds$$

$$= \frac{1}{|\lambda_i|^2} \|f_3\|_{L^2(H,(a_3,+\infty))}^2 < \infty.$$

The above simple calculations are shown that $u_1(\lambda;\cdot) \in L^2(H,(-\infty,a_1))$, $u_3(\lambda;\cdot) \in L^2(H,(a_3,+\infty))$, i.e. $u(\lambda;\cdot) = (u_1(\lambda;\cdot), 0, u_3(\lambda;\cdot)) \in L^2(1,0,1)$ in the case when $\lambda \in \mathbb{C}$, $\lambda_i = Im\lambda < 0$. On the other hand it can be verified that the function $u(\lambda;\cdot)$ satisfy the equation $L_{W_1} u = \lambda u(\lambda;\cdot) + f$ and $u_3(a_3) = W_1 u_1(a_1)$.

Hence the following result has been proved that for the resolvent set $\rho(L_{W_1})$

$$\rho(L_{W_1}) \supset \{\lambda \in \mathbb{C}: Im\lambda \neq 0\}.$$

Now will be researched continuous spectrum $\sigma_c(L_{W_1})$ of the extension $L_{W_1}$.

For the $\lambda \in \mathbb{C}, \lambda_i = Im\lambda > 0$ the norm of resolvent operator $R_\lambda(L_{W_1})$ of the $L_{W_1}$ is in form

$$\|R_\lambda(L_{W_1})f(t)\|_{L^2}^2 = \left\| e^{i(\tilde{A}_1-\lambda)(t-a_1)} f_1^* + i\int_t^{a_1} e^{i(\tilde{A}_1-\lambda)(t-s)} f_1(s) ds \right\|_{L^2(H,(-\infty,a_1))}^2$$

$$+ \left\| i \int_t^{\infty} e^{i(\tilde{A}_3-\lambda)(t-s)} f_3(s) ds \right\|_{L^2(H,(a_3,+\infty))}^2, \quad f = (f_1, 0, f_3) \in L^2(1,0,1).$$

Then it is clear that for any $f = (f_1, 0, f_3) \in L^2(1,0,1)$ is true

$$\|R_\lambda(L_W)f(t)\|_{L^2}^2 \geq \left\| i \int_t^{\infty} e^{i(\tilde{A}_3-\lambda)(t-s)} f_3(s) ds \right\|_{L^2(H,(a_3,+\infty))}^2.$$

The vector functions $f^*(\lambda;t)$ in form $f^*(\lambda;t) = (0,0, e^{i(\tilde{A}_3-\bar\lambda)t} f_3)$, $\lambda \in \mathbb{C}$, $\lambda_i = Im\lambda > 0$, $f_3 \in H$ belong to $L^2(1,0,1)$. Indeed,

$$\|f^*(\lambda;t)\|_{L^2}^2 = \int_{a_3}^{\infty} \|e^{i(\tilde{A}_3-\bar\lambda)t} f_3\|_H^2 dt = \int_{a_3}^{\infty} e^{-2\lambda_i t} dt\, \|f_3\|_H^2 = \frac{1}{2\lambda_i} e^{-2\lambda_i a_3} \|f_3\|_H^2 < \infty.$$

For the such functions $f^*(\lambda;\cdot)$ we have

$$\|R_\lambda(L_{W_1})f^*(\lambda;t)\|_{L^2(H,(a_3,+\infty))}^2 \geq \left\| i\int_t^{\infty} e^{i(\tilde{A}_3-\lambda)(t-s)} e^{i(\tilde{A}_3-\bar\lambda)s} f_3 ds \right\|_{L^2(H,(a_3,+\infty))}^2$$

$$= \left\| \int_t^{\infty} e^{-i\lambda t} e^{-2\lambda_i s} e^{i\tilde{A}_3 t} f_3 ds \right\|_{L^2(H,(a_3,+\infty))}^2 = \left\| e^{-i\lambda t} e^{i\tilde{A}_3 t} \int_t^{\infty} e^{-2\lambda_i s} f_3 ds \right\|_{L^2(H,(a_3,+\infty))}^2$$

$$= \left\| e^{-i\lambda t} \int_t^{\infty} e^{-2\lambda_i s} ds \right\|_{L^2(H,(a_3,+\infty))}^2 \|f_3\|_H^2 = \frac{1}{4\lambda_i^2} \int_{a_3}^{\infty} e^{-2\lambda_i t} dt\, \|f_3\|_H^2 = \frac{1}{8\lambda_i^3} e^{-2\lambda_i a_3} \|f_3\|_H^2.$$

From this

$$\|R_\lambda(L_{W_1})f^*(\lambda;\cdot)\|_{L^2} \geq \frac{e^{-\lambda_i a_3}}{2\sqrt{2}\lambda_i\sqrt{\lambda_i}} \|f\|_H = \frac{1}{2\lambda_i} \|f^*(\lambda;\cdot)\|_{L^2}$$

i.e. for $\lambda_i = Im\lambda > 0$ and $f \neq 0$ are valid

$$\frac{\|R_\lambda(L_{W_1})f^*(\lambda;\cdot)\|_{L^2}}{\|f^*(\lambda;\cdot)\|_{L^2}} \geq \frac{1}{2\lambda_i}.$$

On the other hand it is clear that

$$\|R_\lambda(L_{W_1})\| \geq \frac{\|R_\lambda(L_{W_1})f^*(\lambda;\cdot)\|_{L^2}}{\|f^*(\lambda;\cdot)\|_{L^2}}, f_3 \neq 0.$$

Consequently, we have

$$\|R_\lambda(L_{W_1})\| \geq \frac{1}{2\lambda_i} \text{ for } \lambda \in \mathbb{C}, \lambda_i = Im\lambda > 0.$$

Furthermore, here spectrum of selfadjoint extensions of the minimal operator $L_0(0,1,0)$ will be investigated.

**Theorem 3.3.** The spectrum of the selfadjoint extension $L_{W_2}$ of the minimal operator $L_0(0,1,0)$ in the Hilbert space $L^2(0,1,0)$ is in form

$$\sigma(L_{W_2}) = \left\{\lambda \in \mathbb{R} : \lambda = \frac{1}{b_2-a_2}arg\mu + \frac{2n\pi}{b_2-a_2}, n \in \mathbb{Z}, \mu \in \sigma(W_2^* e^{i\tilde{A}_2(b_2-a_2)}), 0 \leq arg\mu < 2\pi\right\}$$

**Proof.** The general solution of the following problem to spectrum of the selfadjoint extension $L_{W_2}$,

$$\tilde{l}_2(u_2) = iu_2' + \tilde{A}_2 u_2 = \lambda u_2 + f_2, \quad u_2, f_2 \in L^2(H, (a_2, b_2))$$
$$u_2(b_2) = W_2 u_2(a_2), \quad \lambda \in \mathbb{R}$$

is in form

$$u_2(t) = e^{i(\tilde{A}_2-\lambda)(t-a_2)} f_2^* + \int_{a_2}^{t} e^{i(\tilde{A}_2-\lambda)(t-s)} f_2(s) ds, a_2 < t < b_2,$$

$$(e^{i\lambda(b_2-a_2)} - W_2^* e^{i\tilde{A}_2(b_2-a_2)}) f_2^* = W_2^* e^{i\lambda(b_2-a_2)} \int_{a_2}^{b_2} e^{i(\tilde{A}_2-\lambda)(b_2-s)} f_2(s) ds$$

From this it is implies that $\lambda \in \sigma(L_{W_2})$ if and only if $\lambda$ is a solution of the equation $e^{i\lambda(b_2-a_2)} = \mu$, where $\mu \in \sigma(W_2^* e^{i\tilde{A}_2(b_2-a_2)})$. From this it is obtained that
$\lambda = \frac{1}{b_2-a_2}arg\mu + \frac{2n\pi}{b_2-a_2}, n \in \mathbb{Z}, \mu \in \sigma(W_2^* e^{i\tilde{A}_2(b_2-a_2)})$.

**Theorem 3.4.** Spectrum $\sigma(L_{W_1 W_2})$ of any selfadjoint extension $L_{W_1 W_2} = L_{W_1} \oplus L_{W_2}$ cioncides with $\mathbb{R}$.

**Proof.** Validity of this assertion is a simple result of the following claim that a proof of which it is clear. If $S_1$ and $S_2$ two linear closed operators in any Hilbert spaces $H_1$ and $H_2$ respectively, then

$$\sigma_p(S_1 \oplus S_2) = \sigma_p(S_1) \cup \sigma_p(S_2),$$
$$\sigma_c(S_1 \oplus S_2) = (\sigma_p(S_1) \cup \sigma_p(S_2))^c \cap (\sigma_r(S_1) \cup \sigma_r(S_2))^c \cap (\sigma_c(S_1) \cup \sigma_c(S_2))$$

Later on, note that for the singular differential operators for n-th order in scalar case in the another type singularity in the finite interval has been reseached in [9].

**Example:** By the last theorem the spectrum of following boundary value problem
$$i\frac{\partial u(t,x)}{\partial t} + sgnt \frac{\partial^2 u(t,x)}{\partial x^2} = f(t,x), |t| > 1, x \in [0,1],$$
$$i\frac{\partial u(t,x)}{\partial t} - \frac{\partial^2 u(t,x)}{\partial x^2} = f(t,x), \quad |t| < 1/2, x \in [0,1],$$

$$u(1/2,x) = e^{i\psi}u(-1/2,x), \psi \in [0,2\pi)$$
$$u(1,x) = e^{i\varphi}u(-1,x), \varphi \in [0,2\pi),$$
$$u'_x(t,0) = u'_x(t,1) = 0, |t| > 1, |t| < 1/2$$

in the space $L^2\big((-\infty,-1) \times (0,1)\big) \oplus L^2\big((-1/2,1/2) \times (0,1)\big) \oplus L^2\big((1,\infty) \times (0,1)\big)$ coincides with $\mathbb{R}$.